\documentclass[12pt]{article}
\usepackage{amssymb}
\usepackage{amsmath,amsfonts,amsthm}
\begin{document}
%These are the macros needed for the Rings.tex
%-----------------------------------------------------------------------------
%                           New commands
\newcommand{\fab}[2]{\langle {#1}_1,{#1}_2,\ldots ,{#1}_{#2} \rangle}
\newcommand{\abs}[1]{\:|#1|}    \newcommand{\Ker}{{\rm Ker\,}}
\newcommand{\ol}{\overline}     \newcommand{\Z}{ Z\!\!\!Z}
\newcommand{\C}{ l\!\!\!C}       \newcommand{\IH}{ I\!\!H}
\newcommand{\R}{ I\!\!R}         \newcommand{\N}{ I\!\!N}
\newcommand{\Q}{ l\!\!\!Q}       \newcommand{\notsubset}{ /\!\!\!\!\!\!\subset }
\newcommand{\id}{\mbox{id}}     \newcommand{\notdivide}{\!\not |\,}
\newcommand{\pf}{{\bf Proof: }}    \newcommand{\F}{ I\!\!F}
 \newcommand{\dg}{\mbox{deg}\,}
\newcommand{\st}{\; \vline\;}   \newcommand{\norm}[1]{\:\parallel #1 \parallel}
\newcommand{\orb}{\mbox{ \bf Orb}}  \newcommand{\stab}{\mbox{ \bf Stab}}
\newcommand{\rhdn}{\rhd_{_{\!\!\!\!\!\!\neq }}}
%-----------------------------------------------------------------------------
\newcommand{\ba}{$$\begin{array}}\newcommand{\ea}{\end{array}$$}
\newcommand{\bea}{\begin{eqnarray*}}\newcommand{\eea}{\end{eqnarray*}}
\newcommand{\be}{\begin{equation}}\newcommand{\ee}{\end{equation}}
\newcommand{\bd}{\begin{definition}} \newcommand{\ed}{\end{definition}}
\newcommand{\brs}{\begin{remarks}\rm}   \newcommand{\ers}{\end{remarks}}
\newcommand{\br}{\begin{remark}\rm}     \newcommand{\er}{\end{remark}}
\newcommand{\bt}{\begin{theorem}}       \newcommand{\et}{\end{theorem}}
\newcommand{\bl}{\begin{lemma}}         \newcommand{\el}{\end{lemma}}
\newcommand{\bco}{\begin{corollary}}    \newcommand{\eco}{\end{corollary}}
\newcommand{\bp}{\begin{proposition}}   \newcommand{\ep}{\end{proposition}}
\newcommand{\bo}{\begin{observation}\rm}\newcommand{\eo}{\end{observation}}
\newcommand{\bex}{\begin{examples}\rm}   \newcommand{\eex}{\end{examples}}
\newcommand{\bos}{\begin{observations}\rm}\newcommand{\eos}{\end{observations}}
\newcommand{\bx}{\begin{example}\rm}   \newcommand{\ex}{\end{example}}
\newcommand{\bexe}{\begin{exercise}\rm}   \newcommand{\eexe}{\end{exercise}}
\newcommand{\bpf}{\begin{proof}\rm}   \newcommand{\epf}{\end{proof}}
%-----------------------------------------------------------------------------
%                    New environment defintions

\newtheorem{definition}{Definition}[section]
\newtheorem{theorem}[definition]{Theorem}
\newtheorem{lemma}[definition]{Lemma}
\newtheorem{corollary}[definition]{Corollary}
\newtheorem{proposition}[definition]{Proposition}
\newtheorem{remarks}[definition]{Remarks}
\newtheorem{remark}[definition]{Remark}
\newtheorem{observation}[definition]{Observation}
\newtheorem{examples}[definition]{Examples}
\newtheorem{observations}[definition]{Observations}
\newtheorem{example}[definition]{Example}
\newtheorem{exercise}[definition]{Exercise}
%-----------------------------------------------------------------------------

\title{On Kuroda's proof of Hilbert's fourteenth problem in dimensions three and four}
\author{Pramod K. Sharma\\ e--mail:
 pksharma1944@yahoo.com\\
 School Of Mathematics, Vigyan Bhawan, Khandwa Road,\\ INDORE--452
017, INDIA.}
\date{}
\maketitle

                  \section*{Abstract} We generalize [3, Lemma 2.2] and [4, Proposition 2.3]
 and deduce a positive result on Hilbert's fourteenth problem. Further, we
 give a relatively transparent and elementary proof of [3, Theorem
 1.1].
  \section{Introduction}
   Let $K$ be a field of characteristic zero. Let $K[{\bf X}]=K[X_1, \cdots , X_n]$
   be the polynomial ring in $n$ variables over $K$ and let $K({\bf X})=K(X_1, \cdots , X_n)$
   be its field of fractions. The fourteenth problem of Hilbert asks whether $K$-algebra
   $L\cap K[\bf X]$ is finitely generated whenever $L$ is a subfield of $K(\bf X)$
   containing $K$. Zariski [9] showed that $L\cap K(\bf X)$ is a finitely generated
  $K$-algebra if the transcendence degree of $L$ over $K$ is $\leq 2$. Thus the problem
   has an affirmative answer for $n\leq 2$. However, a counter example was found by
   Nagata [7] for $n\geq 32$. Roberts [8] constructed a new counter example for $n=7$.
   Freudenburg [6] gave a counter example for $n=6$, and Daigle and Freudenburg [5]
   gave one for $n=5$. Kuroda [3] and [4] has given counter examples in case $n=4$ and
   $n=3$. Thus Hilbert's fourteenth problem has been solved for all $n$. For $n=4$,
   let $\gamma$ and $\delta_{ij}$ be integers for $1\leq i\leq 3$ and $1\leq j\leq 4$
   such that $\gamma$, $\delta_{ij}\geq 1$ for $1\leq i,j\leq 3$ and $\delta_{i4}\geq 0$.
   Assume \ba{ccc} \pi_1 & = &
   X_4^{\gamma}-X_1^{-\delta_{11}}X_2^{\delta_{12}}X_3^{\delta_{13}}X_4^{\delta_{14}},
     \\ \pi_2 & = &
     X_4^{\gamma}-X_1^{\delta_{21}}X_2^{-\delta_{22}}X_3^{\delta_{23}}X_4^{\delta_{24}},
   \\     \pi_3 & = &
     X_4^{\gamma}-X_1^{\delta_{31}}X_2^{\delta_{32}}X_3^{-\delta_{33}}X_4^{\delta_{34}},
     \ea
     and let $K(\pi)$ be the field of fractions of $ K[\pi] = K[\pi_1,\pi_2,
     \pi_3]$. Then Kuroda [3, Theorem 1.1] proves: If $$
     \frac{\delta_{11}}{\delta_{11}+\min(\delta_{21},\delta_{31})} +
     \frac{\delta_{22}}{\delta_{22}+\min(\delta_{32},\delta_{12})} +
     \frac{\delta_{33}}{\delta_{33}+\min(\delta_{13},\delta_{23})}
      < 1 \,\,\, \cdots(\ast), $$ then $K(\pi)\cap K[\bf X]$ is not a finitely
     generated $K$-algebra. Moreover, $K(\pi)\cap K[\bf X]$ is not
     contained in the kernel $K[\bf X]^D$ for any locally nilpotent
     derivation $D$ of $K[\bf X]$. Further, for $n=3$, let $\gamma$ and
     $\delta_{ij}$ be positive integers for $i,j=1,2$ and let
    \[\begin{array}{ccl}
     \pi_1 & = & X_1^{\delta_{21}}X_2^{-\delta_{22}} - X_1^{-\delta_{11}}X_2^{\delta_{12}} \\
    \pi_2 & = & X_3^{\gamma} - X_1^{-\delta_{11}}X_2^{\delta_{12}} \\
     \pi_3 & = & 2
     X_1^{\delta_{21}-\delta_{11}}X_2^{\delta_{12}-\delta_{22}} -
     X_1^{-2\delta_{11}}X_2^{2\delta_{12}}.
     \end{array}\]
     Then Kuroda [4, Theorem 1.1] states that: if
     \[ \frac{\delta_{11}}{\delta_{11}+\delta_{21}}+
     \frac{\delta_{22}}{\delta_{22}+\delta_{12}} < \frac{1}{2}
     \;\;\;\; \cdots \cdots(\ast\ast),\] then for
     $K(\pi)=K(\pi_1,\pi_2,\pi_3), \;K(\pi)\cap K[X_1,X_2,X_3]$ is not
     a finitely generated $K$-algebra. \\ In proving the above results,
     the first main step is to show that $K(\pi)\cap K[\bf
     X]=K[\pi]\cap k[\bf X]$, and then show that this algebra is not
     finitely generated. In this note we shall show that the statement
     $K(\pi)\cap K[\bf X]=K[\pi]\cap k[\bf X]$ has common genus and
      deduce this fact, i.e., [3, Lemma 2.2] and [4, Proposition 2.3] are
     consequences of one algebraic result under a condition weaker
     than $(*)$ and $(**)$. To be precise: From the proof in[3] it is
      clear that the condition $(*)$ implies that the matrix \[ T_3 =
      \left(\begin{array}{ccc}
      -\delta_{11} & \delta_{12} & \delta_{13} \\
       \delta_{21} & -\delta_{22} & \delta_{23} \\
    \delta_{31} & \delta_{32} & -\delta_{33}
    \end{array}\right)\] is invertible. We also note that the
     condition $(**)$ implies invertibility of the matrix \[T_2 =
    \left(
    \begin{array}{cc}
       -\delta_{11} & \delta_{12} \\
    \delta_{21} & -\delta_{22}
   \end{array}\right)\] The converse is not true in both the cases. We
   prove [3, Lemma 2.2] and [4, Proposition 2.3] under the assumption
   that the corresponding matrices are invertible. We further show
   that the condition in [3, Theorem 1.1] is essential by giving an
   example wherein $K[\pi]\cap K[\bf X]$ is finitely generated if the
   condition $(*)$ fails. We also deduce a positive result [Theorem 2.8]
   on Hilbert's fourteenth problem using our algebraic result and give a
   relatively transparent and easy proof of [3, Theorem 1.1 ]. All rings
   are commutative with identity $(\neq 0)$, and for any ring $R$, $R^*$
  denotes its group of units.
\section{Main Results}
We first record the following result due to Bass. \bl
\label{h-21}[1, Proposition D.1.7]: Let $A\subset B$ be a ring
extension, where $A$ and $B$ are domains with quotient fields
$Q(A)$ and $Q(B)$ respectively. Suppose that $A$ is a unique
factorization domain and that $B^*\cap A\subset A^*$. If $B$ is
flat over $A$, then $Q(A)\cap B=A$. \el
 \bco Let $A$ and $B$ be
integral domains contained in a common integral domain $C$. If $A$
is unique factorization domain and $C$ is flat over $A$ such that
$C^*\cap A\subset A^*$, then $Q(A)\cap B =A\cap B$ where $Q(A)$ is
the field of fractions of $A$. \eco
 \bpf The proof is immediate
from the Lemma \ref{h-21}. \epf
 \br In the Lemma \ref{h-21}, we
can replace the condition $B^*\cap A\subset A^*$ by the assumption
that $B$ is integral over $A$. \er \bl\label{2.4} Let $A\subset B$
be integral domains, where $A$ is a unique factorization domain.
If $B$ is faithfully flat over $A$, then $Q(A)\cap B=A$ where
$Q(A)$ is the field of fractions of $A$. \el \bpf The proof will
 follow
 from  Lemma \ref{h-21}, if $B^*\cap A\subset A^*$. Let $\lambda
\in B^*\cap A$. As $A$ is an integral domain,  \be \label{h-22}
0\longrightarrow A \stackrel{ m_\lambda} {\longrightarrow} A
\stackrel{\eta}{\longrightarrow}A/\lambda A \longrightarrow 0 \ee
is an exact sequence of $A$-modules where $m_\lambda$ is
multiplication by $\lambda$, and $\eta$ is the quotient map. Since
$B$ is faithfully flat, from the exact sequence \ref{h-22}, we get
the exact sequence $$ 0\longrightarrow B \stackrel{ m_\lambda}
{\longrightarrow} B \stackrel{\eta}{\longrightarrow}B /\lambda B
\longrightarrow 0 $$ As $\lambda \in B^*$, $B/\lambda B=0$. Thus
$A/\lambda A\otimes B\cong  B/\lambda B=(0)$. Consequently, as $B$
is faithfully flat, $A/\lambda A=(0)$. Thus $\lambda \in A^*$, and
the result follows. \epf \bt\label{2.5} Let $K$ be a field and
$X_i;$ $i=1,2,\cdots , n$ be algebraically independent over $K$.
Let $\pi_i=X_n^{\gamma_i}-X_1^{\alpha_{i1}}X_2^{\alpha_{i2}}\cdots
X_{i-1}^{\alpha_{ii-1}}X_i^{-\alpha_{ii}}\cdots X_n^{\alpha_{in}}$
for $i=1,2,\cdots , n-1 $ where $\gamma_i\geq 1$, $\alpha_{ij}\geq
1$ for all $1\leq i, j\leq n-1$ and $\alpha_{in}\geq 0$. If for
the matrix \[ T_{n-1} = \left(\begin{array}{cccc}
 -\alpha_{11} & \alpha_{12} & \cdots & \alpha_{1n-1} \\
  \alpha_{21} & -\alpha_{22} & \cdots & \alpha_{2n-1} \\ \cdots & \cdots
  & \cdots & \cdots  \\  \alpha_{(n-1)1} & \cdots & \cdots & -\alpha_{(n-1)(n-1)}
\end{array}\right),\] $\det T_{n-1}\neq 0$, then \\ $(i)$ $\pi_1, \pi_2,
\cdots \pi_{n-1}, X_n $ are algebraically independent over $K$. \\
$(ii)$ The ring extension $$K[\pi]=K[\pi_1,\cdots ,
\pi_{n-1}]\subset K[X^{\pm 1}_1, \cdots , X_n^{\pm 1}] $$ is flat.
\et \bpf $(i)$ We have $$K(\pi_1,\cdots , \pi_{n-1})\subset K(X_n,
\pi_1,\cdots , \pi_{n-1})=K(X_n, M_1, \cdots M_{n-1})=L\mbox{
(say) } $$ where $M_i=X_1^{\alpha_{i1}}\cdots
X_{i-1}^{\alpha_{ii-1}} X_{i}^{-\alpha_{ii}}\cdots
X_{n-1}^{\alpha_{in-1}}$ for $ i = 1,2,\cdots ,n-1$. We shall show
that $K(X_1, \cdots , X_n)\,|\,L$ is an algebraic field extension.
Since $\det T_{n-1}\neq 0$, for any $i=1,2, \cdots , n-1$ there
exists a row vector $f_i\in \Q^{n-1}$ such that \ba{rl} &
f_iT_{n-1}=e_i=(0, \cdots, 0, 1^{i^{th}}, 0, \cdots , 0) \\
\Rightarrow & (mf_i)T_{n-1}=me_i\,\mbox{ for all } m\in \Z. \ea
Choose $m_i>0$ such that $m_if_i\in \Z^{n-1}$. If $m_if_i=(s_1,
\cdots , s_{n-1})\in \Z^{n-1}$, then $$M_1^{s_1}\cdots
M_{n-1}^{s_{n-1}}=X_i^{m_i}.$$ Therefore $X_i^{m_i}\in L$. Thus
$K(X_1, \cdots , X_n)$ is algebraic over $L$. This proves $(i)$.

$(ii)$ From $(i)$, $X_n,\pi_1,\cdots,\pi_{n-1}$ are algebraically
independent over $K$. Therefore $K[\pi]\subset
K[\pi,X_n,X_n^{-1}]$ is flat. Also
$$K[\pi,X_n,X_n^{-1}]=K[X_n,X_n^{-1},M_1,\cdots,M_{n-1}] \subset
K[X_n,X_n^{-1},M_1^{\pm 1},\cdots,M_n^{\pm1}]=B$$ is flat. Now,
note that the Laurent polynomial ring $K[X_1^{\pm
1},\cdots,X_n^{\pm 1}]=R$ is a graded ring over the free abelian
group $(\Z^n,+)$ where the homogeneous component corresponding to
$(\alpha_1,\cdots,\alpha_n)\in\Z^n$ is $KX_1^{\alpha_1}\cdots
X_n^{\alpha_n}$. Then the ring $B$ is a graded subring of $R$ over
the subgroup $H$ of $\Z^n$ generated by the set \[\begin{array}{l}
  \{(-\alpha_{11},\alpha_{12},\cdots,\alpha_{1(n-1)},0),
  (\alpha_{21},-\alpha_{22},\cdots,\alpha_{2(n-1)},0),\cdots, \\
  (\alpha_{(n-1)1},\cdots,\alpha_{(n-1)(n-2)},-\alpha_{(n-1)(n-1)},0),
  (0,\cdots,0,1)\}
\end{array}\] Therefore $R$ is a free $B$-module. Hence $(ii)$ follows. \epf
\bp\label{2.6} Let $K$ be a field and $X_i;\;i=1,2,\cdots,n$ be
algebraically independent over $K$. Let \[\pi_i=X_n^{\gamma_i}-
X_1^{\alpha_{i1}}X_2^{\alpha_{i2}}\cdots
X_{(i-1)}^{\alpha_{i(i-1)}}X_i^{-\alpha_{ii}}\cdots
X_n^{\alpha_{in}}\] for $i=1,2,\cdots,n-1$ where $\gamma_i\geq 1$,
$\alpha_{ij}\geq 1$ for all $1\leq i,j\leq n-1$ and
$\alpha_{in}\geq 0.$ If for the matrix \[T_{n-1}=\left(\begin{array}{cccc}
  -\alpha_{11} & \alpha_{12} & \cdots & \alpha_{1(n-1)} \\
  \alpha_{21} & -\alpha_{22} & \cdots & \alpha_{2(n-1)} \\
  \cdots & \cdots & \cdots & \cdots \\
  \alpha_{(n-1)1} & \cdots & \cdots & -\alpha_{(n-1)(n-1)}
\end{array}\right)\] $\det T_{n-1} \neq 0,$ then for $K[\pi]=K[\pi_1,\cdots,\pi_{n-1}]$
\[K(\pi)\cap K[X_1^{\pm 1},\cdots,X_n^{\pm 1}]=K[\pi]\] where $K(\pi)$ is the field of
fractions of $K[\pi]$. \ep \bpf By Theorem \ref{2.5}, the ring
extention $K[\pi]\subset K[X_1^{\pm 1},\cdots, X_n^{\pm 1}]$ is
flat, and $\pi_1,\cdots,\pi_{n-1},X_n$ are algebraically
independent over $K$. As $\pi_1,\cdots,\pi_{n-1}$ are
algebraically independent over $K$, $K[\pi_1,\cdots,\pi_{n-1}]$ is
a unique factorization domain. As in the proof of the Theorem
\ref{2.5}$(ii)$, it is clear that $R=K[X_1^{\pm 1},\cdots,
X_n^{\pm 1}]$ is free over its subring $B=K[X_n^{\pm 1},M_1^{\pm
1},\cdots,M_{n-1}^{\pm 1}].$ Therefore $R$ is faithfully flat over
$B$. By Lemma \ref{2.4}, $R^*\cap B\subset B^*$. Hence
\[\begin{array}{rrr}
  & R^*\cap B\cap K[\pi]\;\subset & B^*\cap K[\pi] \\
  \Longrightarrow & R^*\cap K[\pi]\;\subset & B^*\cap K[\pi]
\end{array}\] since $K[\pi]\subset B.$ Now the result will follow from
Lemma \ref{h-21} if we show that $B^*\cap K[\pi]\subset
K[\pi]^*=K^*.$ Since $X_n,\pi_1,\cdots,\pi_{n-1}$ are
algebraically independent over $K$,
\[K[\pi_1,\cdots,\pi_{n-1},X_n^{\pm 1}]=K[M_1,\cdots,M_{n-1},X_n^{\pm
1}]\] is a unique factorization domain. We also note that
$M_1,\cdots,M_{n-1},X_n$ are algebraically independent over $K$.
Hence $B=K[M_1,\cdots,M_{n-1},X_n,1/f],$ where
$f=M_1.M_2.\cdots.M_{n-1}.X_n$, is a unique factorization domain.
Let for $\lambda\in K[\pi]$, $\lambda$ is unit in $B$. Then there
exists $\mu\in K[M_1,\cdots,M_{n-1},X_n]$ such that
$\lambda\mu=f^m$ for some $m\geq 0$. Note that $K[\pi]\subset
K[M_1,\cdots,M_{n-1},X_n]$. Therefore $\lambda$ divides $f^m$ in
$K[M_1,\cdots,M_{n-1},X_n]$. Consequently \[\lambda =
aM_1^{\alpha_1}\cdots M_{n-1}^{\alpha_{n-1}}X_n^{\alpha_n}
\;\;\;(\alpha_i\geq 0)\] for some $a(\neq 0)\in K$. We have $\pi_i
= X_n^{\gamma_i}-M_iX_n^{\alpha_{in}}.$ Therefore $\lambda =
aM_1^{\alpha_1}\cdots M_{n-1}^{\alpha_{n-1}}X_n^{\alpha_n} \in
K[\pi]$ implies \[\alpha_n =0 =\alpha_1 = \cdots = \alpha_{n-1}.\]
This can be seen by putting $M_i=X_n^{p^i}$ for a sufficiently
large prime $p$ in the equation \[\lambda = aM_1^{\alpha_1}\cdots
M_{n-1}^{\alpha_{n-1}}X_n^{\alpha_n}.\] Hence $\lambda\in K^*$ and
the result follows.  \epf \br\label{2.7} From the above result,
[3, Lemma 2.2] is immediate.\er We now prove a positive result on
Hilbert's fourteenth problem using our results. \bt\label{2.8} Let
$K$ be a field and let $X_1, X_2,\cdots,X_n$ be algebraically
independent over $K$. Let $M_1,M_2,\cdots,M_t$ be monomials in
$X_1^{\pm 1},\cdots,X_n^{\pm 1}$ algebraically independent over
$K$. Then
\[K(M_1,\cdots,M_t)\cap K[X_1,\cdots,X_n]\] is a finitely
generated $K$-algebra.\et \bpf Since $M_1,M_2,\cdots,M_t$ are
algebraically independent over $K$,\[K[M_1^{\pm 1}, M_2^{\pm
1},\dots,M_t^{\pm 1}]= K[M_1,M_2,\cdots,M_t,1/f],\] where $f=
M_1M_2\cdots M_t$, is a unique factorization domain. Further, \[
K[ M_1^{\pm 1}, M_2^{\pm 1},\cdots,M_t^{\pm 1}]\subset K[X_1^{\pm
1}, \cdots,X_n^{\pm 1}]\] is a faithfully flat ring extension (in
fact a free extension). Hence, by Lemma \ref{2.4},
\[K(M_1,M_2,\cdots ,M_t)\cap K[X_1^{\pm 1},\cdots,X_n^{\pm 1}] =
K[M_1^{\pm 1}, M_2^{\pm 1},\cdots,M_t^{\pm 1}].\] Therefore
\[K(M_1,M_2,\cdots ,M_t)\cap K[X_1,\cdots,X_n] = K[M_1^{\pm 1},
M_2^{\pm 1},\cdots,M_t^{\pm 1}]\cap K[X_1,\cdots,X_n].\] Now, as
$M_1,M_2,\cdots,M_t$ are algebraically independent over $K$, an
element \[f=\sum\lambda_{\alpha_1,\cdots,\alpha_t} M_1^{\alpha_1}
\cdots M_t^{\alpha_t}\;\;\;(\lambda_{\alpha_1,\cdots,\alpha_t} \in
K,\; \alpha_i\in\Z)\;\cdots(***)\] is in $K[M_1^{\pm 1},\cdots,
M_t^{\pm 1}]\cap K[X_1,\cdots,X_n]$ if and only if $M_1^{\alpha_1}
M_2^{\alpha_2}\cdots M_t^{\alpha_t}$ is a monomial in
$X_1,\cdots,X_n$. Let $M_i=X_1^{\gamma_{i1}}\cdots
X_n^{\gamma_{in}}$ for $i=1,\cdots,t$ and let
\[U=\left(\begin{array}{ccc}
  \gamma_{11} & \cdots & \gamma_{1n} \\
  \gamma_{21} & \cdots & \gamma_{2n} \\
  \cdots & \cdots & \cdots \\
  \gamma_{n1} & \cdots & \gamma_{nn}
\end{array}\right).\] Thus in the equation $(***)$, all power-tuples $(\alpha_1,\cdots,
\alpha_t)$ are in the
semi-group\[S=\{(\beta_1,\cdots,\beta_t)\in\Z^t\st(\beta_1,\cdots,\beta_t
)U=(a_1,\cdots,a_t) \mbox{ where }a_i\geq 0\}.\] By Gordan's Lemma
[2], $S$ is finitely generated. Hence the result follows. \epf

  We shall now deduce [4, Proposition 2.3]. Observe that the condition $(\ast\ast)$,
  i.e.,
      \[
\frac{\delta_{11}}{\delta_{11}+\delta_{21}}+
\frac{\delta_{22}}{\delta_{22}+\delta_{12}} < \frac{1}{2} \]
implies  the determinant of the matrix
 \[ T_2 = \left(\begin{array}{cc}
 -\delta_{11} & \delta_{12}  \\
  \delta_{21} & -\delta_{22} \\

\end{array}\right)\] is non-zero. If not, then there exist $(t,s)\neq
(0,0)\in\Z^2$ such that $(t,s)T_2=0.$ Hence
$-t\delta_{11}+s\delta_{21}=0$ and $t\delta_{12}-s\delta_{22}=0.$
It is easy to see that $t\geq 0$ if and only if $s\geq 0$.
Moreover, $t=0$ if and only if $s=0$. Let $t\geq s>0$. Then
$\frac{s}{t} \leq 1$. Hence
$\delta_{12}=\frac{s}{t}\,\delta_{22}\leq \delta_{22}.$ This gives
\[\frac{\delta_{22}}{\delta_{22}+\delta_{12}}\geq
\frac{\delta_{22}}{\delta_{22}+\delta_{22}}= \frac{1}{2}.\]
Further, if $s\geq t>0,$ then
\[\frac{\delta_{11}}{\delta_{11}+\delta_{21}} \geq \frac{\delta_{11}}{\delta_{11}
+\delta_{11}} =\frac{1}{2}.\] Thus in any case we arrive at a
contradiction to our assumption. Hence $\det T_2\neq 0.$ The
converse is not true follows by taking $\delta_{11}=3,
\delta_{12}= \delta_{22}=\delta_{21}=1.$ Now assume $\det T_2\neq
0.$ To prove [4, Proposition 2.3] under this assumption we first
note the following: \bl\label{2.13} Let $R$ be a unique
factorization domain and $K$ its field of fractions. Let $L|K$ be
an algebraic field extension. If $x\in L$ is integral over $R$,
then $R[x]$ is finitely generated free $R$-module.\el \bpf Let
$f(Z)\in R[Z]$ be a monic polynomial of least degree such that
$f(x)=0.$ Then $f(Z)$ is irreducible in $R[Z]$ and has content 1.
Therefore $f(Z)$ is irreducible in $K[Z]$. Since $f(Z)$ is monic,
it is easy to see that $f(Z)K[Z]\cap R[Z]=f(Z)R[Z].$ Therefore
$R[Z]/(f(Z))\cong R[x]$ is a finitely generated free $R$-module.
\epf \bt\label{2.14} Let $X_1,X_2,X_3$ be algebraically
independent over a field $K$. Let
\[\begin{array}{ccl}
  \pi_1 & = & X_1^{\delta_{21}}X_2^{-\delta_{22}}- X_1^{-\delta_{11}} X_2^{\delta_{12}}, \\
  \pi_2 & = & X_3^{\gamma}-X_1^{-\delta_{11}}X_2^{\delta_{12}}, \\
  \pi_3 & = &
  2X_1^{\delta_{21}-\delta_{11}}X_2^{\delta_{12}-\delta{22}} -
  X_1^{-2\delta_{11}}X_2^{2\delta_{12}},
\end{array}\] where $\gamma$ and $\delta_{ij}$ are positive
integers for $i,j=1,2.$ Then if for \[T_2 = \left(
\begin{array}{cc}
  -\delta_{11} & \delta_{12} \\
  \delta_{21} & -\delta_{22}
\end{array}\right)\] $\det T_2\neq 0,\;\pi_1,\pi_2,\pi_3$ are
algebraically independent over $K$. Moreover,
$K(\pi_1,\pi_2\pi_3)\cap K[X_1,X_2,X_3]=K[\pi_1,\pi_2,\pi_3] \cap
K[X_1,X_2,X_3].$\et \bpf Since $\det T_2\neq 0$, there exists a
row vector $(m,n)\in \Z^2$ such that
\[(m,n)\left(\begin{array}{cc}
  -\delta_{11} & \delta_{12} \\
  \delta_{21} & -\delta_{22}
\end{array}\right)=(k,0)\] for some $k>0$. Therefore \[ (X_1^{-\delta_{11}} X_2^{
\delta_{12}})^m(X_1^{\delta_{21}}X_2^{-\delta_{22}})^n=X_1^k\]
Thus \[X_1^k\in K[\pi_1,\pi_2,X_3]= K[X_1^{-\delta_{11}} X_2^{
\delta_{12}}, X_1^{\delta_{21}}X_2^{-\delta_{22}}, X_3].\]
Similarly, there exists an integral $l>0$ such that $X_2^l\in
K[\pi_1,\pi_2,X_3].$ Consequently $K[X_1,X_2,X_3]$ is integral
over $K[\pi_1,\pi_2,X_3].$ Thus \[Tr.deg._KK(\pi_1,\pi_2,X_3) =3
\;\;\;\;\cdots (i)\] Hence \[2\leq Tr.deg._KK(\pi_1,\pi_2,\pi_3)
\leq 3,\;\;\;\;\cdots(ii)\] since $K[\pi_1,\pi_2,\pi_3]\subset
K[\pi_1,\pi_2,X_3].$ Now, note that $\pi_2-\pi_1+X_1^{\delta_{21}}
X_2^{-\delta_{22}}=X_3^{\gamma}$ and
$\pi_1^2+\pi_3=X_1^{2\delta_{21}}X_2^{-2\delta_{22}}.$ Therefore
$(X_3^{\gamma}+(\pi_1-\pi_2))^2 = \pi_1^2+\pi_3.$ Hence,
$X_3^{2\gamma}+2(\pi_1-\pi_2)X_3^{\gamma}+
(\pi_1-\pi_2)^2-(\pi_1^2+\pi_3) =0$, i.e., $X_3^{\gamma}$ is
integral over $K[\pi_1,\pi_2,\pi_3].$ Therefore, in view of $(i)$
and $(ii)$, $Tr.deg._KK(\pi_1,\pi_2,\pi_3)=3.$ This proves that
$\pi_1,\pi_2,\pi_3$ are algebraically independent over $K$. For
the last part of the statement, note that $K[\pi_1,\pi_2,\pi_3]$
is a unique factorization domain. As $X_3$ is integral over it,
$K[\pi_1,\pi_2,\pi_3,X_3]$ is a free $K[\pi_1,\pi_2,\pi_3]$-module
by Lemma \ref{2.13}. We have $$ K[\pi_1,\pi_2,\pi_3,X_3]= K
[X_1^{-\delta_{11}}X_2^{\delta_{12}},X_1^{\delta_{21}}X_2^{-\delta_{22}},
X_3 ] \subset K[X_1^\pm 1,X_2^\pm 1,X_3^\pm 1].$$  Therefore, as
in the proof of Theorem \ref{2.5}, $K[X_1^{\pm 1}, X_2^{\pm 1},
X_3^{\pm 1}]$ is $K[\pi_1,\pi_2,\pi_3,X_3]$-free. Consequently
$K[X_1^{\pm 1}, X_2^{\pm 1}, X_3^{\pm 1}]$ is
$K[\pi_1,\pi_2,\pi_3]$-free module. Thus by Lemma \ref{2.4}, \[
K(\pi_1,\pi_2,\pi_3) \cap K[X_1,X_2,X_3] =
K[\pi_1,\pi_2,\pi_3]\cap K[X_1,X_2,X_3]\]
 \epf At the end, we show that [3, Theorem 1.1] is not true if instead
 of $(\ast)$ condition we assume that  determinant of
 the matrix $T_{n-1}$ is non-zero. This is  consequence of the
 following Lemma. \bl\label{2.15} Let $K$ be a field of
 characteristic $\neq 2.$ If $X_i\,; i=1,2,3,4$ are  algebraically
 independent over $K$, then \[K[X_1,X_2,X_3,X_4]\cap
 K[X_4-X_1^{-1}X_2X_3, X_4-X_1X_2^{-1}X_3,
 X_4-X_1X_2X_3^{-1}]=K.\]\el \bpf Let $a,b,c,x$ be algebraically independent
 over $K$. Replacing $X_1,X_2,X_3,X_4$ by $ab,bc,ca$ and $x$ respectively,
  the statement of Lemma amounts to proving: \[K[ab,bc,ca,x] \cap K[
  x-c^2,x-a^2,x-b^2] =K.\] We shall prove this result in steps. \\
  Step 1. A monomial $a^ib^jc^k\in K[ab,bc,ca]$ if and only if $i+j+k$
   is an even integer and $i+j\geq k,\,j+k\geq i,\,i+k\geq j.$ \\
   Let us note that the monomials $ab,bc,ca$ satisfy the conditions in the statement.
   Therefore any monomial $a^ib^jc^k$ in $K[ab,bc,ca]$ also satisfies the same
   conditions. Now, let $a^ib^jc^k$ be a non-constant monomial such that $i+j+k$
   is even and $i+j\geq k,\,j+k\geq i,\,i+k\geq j$. We shall show that $a^ib^jc^k
   \in K[ab,bc,ca]$ by induction on $i+j+k$. Clearly $i+j+k \geq 2$. If $i+j+k = 2$,
    then the monomial is either $ ab $ or $ bc $ or $ ca $ and hence is in
    $ K[ab,bc,ca] $. Assume $i+j+k > 2$, and let $i\geq j\geq k $. If $k=0$,
    then $i=j$ and the result follows . Thus assume $ k\geq 1$.
    As  $i+j+k$ is even we can not have $ i=j=k=1 $. Therefore $
    i+j\geq k+2 $. Then $a^ib^jc^k = (a^{i-1}b^{j-1}c^k)(ab) $,
    and $a^{i-1}b^{j-1}c^k $ satisfies the condition of the
    statement . Therefore, by induction, $a^{i-1}b^{j-1}c^k \in
    K[ab, bc, ca]$. Hence the statement is true.\\
    Step 2. $K[ab,bc,ca]\cap K[x-a^2,x-b^2,x-c^2]=K.$\\ We have
    $K[x-a^2,x-b^2,x-c^2]=K[b^2-a^2,c^2-a^2][x-a^2].$ Let
    $f(ab,bc,ca)\in K[ab,bc,ca]$ be a non-constant polynomial in
    $K[b^2-a^2,c^2-a^2][x-a^2]$. Then it is immediate that \[
    f(ab,bc,ca)=g(b^2-a^2,c^2-a^2)\in
    K[b^2-a^2,c^2-a^2]. \;\;\;\cdots(i)\] We can assume that $f$ and
    $g$ are homogeneous of same degree $d$. Let
    $(b^2-a^2)^i(c^2-a^2)^j$ be the term in $g(b^2-a^2,c^2-a^2)$
    with largest $i$. We have $i+j=d$. If $i>d/2$, then by Step 1,
    the monomial $b^{2i}c^{2j}$ is not in $K[ab,bc,ca]$, since
    $2j<2i.$ Therefore $(i)$ does not hold in this case. Similarly
    if $(b^2-a^2)^i(c^2-a^2)^j$ is the term in
    $g(b^2-a62,c^2-a^2)$ with largest $j$ and $j>d/2,\; (i)$ can
    not hold. Thus assume $i=j=d/2$. Then $(b^2-a^2)^{d/2}
    (c^2-a^2)^{d/2}\in K[ab,bc,ca]$. This gives that $a^{2d}\in
    K[ab,bc,ca]$, which is not true by Step 1. Thus $(i)$ does not
    hold in this case as well. Consequently Step 2 is proved. \\
    Step 3. If $Char.K=0$, then the result of the Lemma holds. \\
    We have $K[ab,bc,ca,x]\cap K[x-a^2,x-b^2,x-c^2] =
    K[ab,bc,ca][x]\cap K[b^2-a^2,c^2-a^2][x-a^2]$. Let \[ f
    =\sum_{i=0}^df_ix^i = \sum_{i=0}^dg_i(x-a^2)^i,\] where $f_i
    \in K[ab,bc,ca]$ and $g_i\in K[b^2-a^2,c^2-a^2]$ be a
    non-constant polynomial of least degree $d$ in $x$ in the
    intersection. We shall first see that $d=0.$ Note that the
    rings $K[ab,bc,ca,x]$ and $K[x-a^2,b^2-a^2,c^2-a^2]$ are both
    closed with respect to partial derivative $\delta/\delta x$.
    Thus if $d\geq 1$, then operating $\delta/\delta x$ on the
    both sides we arrive at a contradiction to minimality of $d$.
    Therefore $d=0$, and $f_0=g_0\in K[ab,bc,ca]\cap
    K[b^2-a^2,c^2-a^2]$. Now, the result follows from Step 2. \\
    Step 4. If characteristic of $K$ is $p>3,$ then
    $K[ab,bc,ca][x]\cap K[b^2-a^2,c^2-a^2][x-a^2]=K$. \\ As in
    step 3, let \[f=\sum_{i=0}^d f_ix^i = \sum_{i=0}^dg_i
    (x-a^2)^i\;\;\;\cdots (ii)\] be a non-constant polynomial of
    least degree  d in $x$ in the intersection. If $d=0$, then the
    statement follows by Step 2. Next, if $d=1$ then comparing
    the coefficients of powers of $x$ on both sides we get
    $f_1=g_1$ and $f_0=g_0-g_1a^2.$ By Step 2, $f_1=g_1=\lambda(
    \neq 0)\in K.$ Hence $f_0=g_0-\lambda a^2.$ This implies that
    either $a^2$ or $b^2$ or $c^2$ is in $K[ab,bc,ca].$ This however
    is not true by Step 1. Now assume $d>1.$ Unless all non-zero
    terms other than that of degree 1 have exponent in $x$ divisible by
    $p$, we can argue as in case of $Char.K=0$ and conclude $d=0$.
    In this case the assertion is already proved. Therefore assume
    that all non-zero terms other than those of degree 1 have
    power of $x$ divisible by $p$. Now, let $q=p^e(e\geq 1)$ be the
    highest power of $p$ which  divides all exponents of
    $x$ other than 1 in
    non-zero term. Let $d=mq$. Then \[\begin{array}{ccl}
      f & = & f_0+f_1x+f_qx^q+\cdots +f_dx^{mq} \\
       & = & g_0+g_1(x-a^2)+\cdots +g_d(x-a^2)^{mq} \\
      \Rightarrow &  & f_0=g_0-g_1a^2+\mbox{ terms with higher power of $a^2$, and
      } f_1=g_1. \\

\end{array}\] If $f_1\neq 0$, then by Step 2,
$f_1=g_1=\lambda(\neq 0)\in K.$ Therefore $f_0=g_0-\lambda a^2+$
terms with higher power of $a^2$. This, however, is not true as
seen above. Thus $f_1=0=g_1.$ Put $Y=X^q.$ Then \[f =f_0+ f_qY+
\cdots +f_dY^m = g_0+g_q(Y-a^{2q}) + \cdots +g_d(Y-a^{2q})^m
\;\;\;\cdots (iii)\] Here $f\in K[Y,ab,ac,bc]\cap K[b^2-a^2,
c^2-a^2][Y-a^{2q}].$ Differentiating both sides of equation
$(iii)$ with respect to $Y$ we get, as in Step 3, $m=0$ since if
$m=1$, then $f_q=g_q,\;f_0=g_0-g_qa^2.$ As $f_q =g_q \in K[
ab,bc,ca]\cap K[b^2-a^2,c^2-a^2]=K.$ By Step 2 $f_q=g_q=\lambda
(\neq 0)\in K$. Then $f_0=g_0-\lambda a^2\in K[ab,bc,ca].$
Therefore, as seen above, $a^2$ or $b^2$ or $c^2$ belong to
$K[ab,bc,ca],$ which is not true in view of Step 1. Thus $m=0$,
and assertion holds by Step 2.\epf \br\label{2.16} The result is
not true if characteristic of $K$ is 2. In this case $(b^2-a^2)
(c^2-a^2)+(x-a^2)^2=(bc)^2-(ac)^2-(ab)^2+x^2 \neq 0$ is in $K[x-
a^2,x-b^2,x-c^2]\cap K[x,ab,bc,ca].$\er
\section{$K[\pi]\cap K[\bf X]$ is not finitely generated}
 In the notations of introduction the main assertion of [ 3,
 Theorem 1.1 ] is that $K[\pi_1,\pi_2,\pi_3]\cap K
 [X_1,X_2,X_3,X_4]$ is not a finitely generated $K-$ algebra. We
 shall give a relatively transparent and easy proof of this
 fact. Put\\
 $Y_1 = X_1^{-\delta_{11}}X_2^{\delta_{12}}X_3^{\delta_{13}}
 X_4^{\delta_{14}},Y_2 =
 X_1^{\delta_{21}}X_2^{-\delta_{22}}X_3^{\delta_{23}}
 X_4^{\delta_{24}},Y_3 =
 X_1^{\delta_{31}}X_2^{\delta_{32}}X_3^{-\delta_{33}}
 X_4^{\delta_{34}} $ and $Y_4 = X_4^{\gamma} $.\\
  Then $K[ \pi ] = K[\pi_1,\pi_2,\pi_3] = K[ Y_4 - Y_1, Y_4 - Y_2,
  Y_4-Y_3]$ and $Y_1, Y_2, Y_3, Y_4 $ are algebraically
  independent by Theorem 2.5(i). As in [ 3 ], let $E$ be the
  $K-$derivation $ E= \Sigma_{i=1}^4 \delta/\delta Y_i $ on $ K[Y]
  = K [Y_1, Y_2, Y_3, Y_4]$. We shall use the elementary facts that
  $ K[Y]^E = K[ Y_4 - Y_1, Y_4 - Y_2, Y_4-Y_3]$ and $ K[Y_2 ,Y_3
  ]^E = K[Y_2 - Y_3]$. Instead of giving complete proofs again,
  we shall prove the main assertions in which our proof differs to
  that in [ 3 ].\\
  The assertion in [ 3, Lemma 2.4 ] is most important for the
  proof of [ 3, Theorem 1.1 ]. Further [3, Lemma 3.1 ] is
  crucial for the proof of  [ 3, Lemma 2.4 ]. We, instead, use the
   following  simple result.
    \bl There exist $ +ve $ integers $p_i ; i=1,2,3 $ such that for
     $ p = p_1 + p_2 + p_3 $, $ p_i \geq p {\xi}_i $ . Moreover, $ f_0 =
    (Y_3-Y_2)^{p_1}(Y_3-Y_1)^{p_2}(Y_2-Y_1)^{p_3} \in K[\bf X] $.
        \el
 \bpf  Since ${\xi}_1 + {\xi}_2 + {\xi}_3 < 1$, there
  exists a $ +ve $ integer $ p $ such that
  $ p(1 - {\xi}_1 - {\xi}_2 - {\xi}_3) \geq 3 $, i.e., $ p -
  p {\xi}_i - p {\xi}_2 - p {\xi}_3 \geq 3 $. Therefore there
  exist $+ve $ integers $p _i ; i = 1,2,3 $ such that
  $ p = p_1 + p_2 + p_3 $ where $ p_i \geq p {\xi}_i $. We have
  $$f_0 =  (Y_3-Y_2)^{p_1}(Y_3-Y_1)^{p_2}(Y_2-Y_1)^{p_3}$$
        $$=  Y_1^{ p_2 + p_3}Y_2^{p_1}(Y_2^{-1}Y_3 -
                1)^{p_1}(Y_1^{-1}Y_3 - 1)^{p_2}(Y_1^{-1}Y_2
                -1)^{p_3}.$$ From this, one can easily see that the
  power of $X_1 $ in any term of $ f_0\in K[ Y_1,Y_2,Y_3 ] $ is
  larger  than $ -(p_2 + p_3)\delta_{11} + p_1  \delta_{21} + (p_1
  - i)\delta_{31} -(p_1 - i)\delta_{21}$ for some $ o \leq i \leq p_1
  $, since clearly the  power of $ X_1$ in  each term of $(Y_1^{-1}Y_3 -
  1)^{p_2}$ and $ (Y_1^{-1}Y_2 - 1)^{p_3}$ is positive. Now
 \[\begin{array}{lcl}
   -(p_2 + p_3)\delta_{11} + p_1  \delta_{21} + (p_1
  - i)\delta_{31} -(p_1 - i)\delta_{21}& = & -(p_2 + p_3)\delta_{11} + (p_1
  - i)\delta_{31} + i\delta_{21}\\ &  \geq & -(p_2 + p_3)\delta_{11} +
  p_1 min(\delta_{21},\delta_{31})\\ & = & -p\delta_{11} + p_1(\delta_{11} +
  min(\delta_{21},\delta_{31}))\\ & = & (p_1 - p{\xi}_1)(\delta_{11} +
  min(\delta_{21},\delta_{31})) \\ &  \geq & 0 \\
   \end{array}\]
  since $p_1 \geq p{\xi}_1$.
    Similarly we can see that powers of $X_2 $ and $X_3$ in $f_0$
  are also $\geq 0$. Thus the result follows.
  \epf
  \br It is easy to prove [3, Lemma 2.4(i)] from the above
  result , but that is not pertinent to the main result.
        \er

   Before we proceed further, to make this proof self contained ,
   we give an alternative proof of [3, Lemma 2.3].
   \bl\label{3.2} Let $ f\in K[ \pi ]\cap K [\bf X ]$ . If $(a_1,a_2,a_3,a_4)\in
   Supp.(f) $ and $a_4 > 0 $, then $ a_1 + a_2 + a_3 > 0 $.
   \el
 \bpf  We have $$ f = \sum{ \lambda_{a,b,c} (Y_4 -Y_1)^a(Y_4 - Y_2)^b(Y_4 -
 Y_3)^c}$$
 where $\lambda_{a,b,c}\in K \mbox{ and } (a,b,c)\in \Z_+^3$. If
 $a + b + c > 0 $, then in a term of the type $ \lambda (Y_4 -Y_1)^a(Y_4 - Y_2)^b(Y_4
 -Y_3)^c $, when expanded, $ \lambda Y_4^{a + b + c} = \lambda
 X_4^{(a + b + c)\gamma} $ occurs, and also $  \lambda
 X_4^{(a + b + c - 1)\gamma}(-Y_1) $ occurs wherein power of
 $X_1$ is $< 0 $. As $ f(X)\in K[X] $, the term $ \lambda
 X_4^{(a + b + c - 1)\gamma}(-Y_1) $ has to cancel out. Thus as $
 Y_1, Y_2, Y_3, Y_4 $ are algebraically independent over $ K $, $ -\lambda
 X_4^{(a + b + c - 1)\gamma}(-Y_1) $ shall also appear in the
 expansion. Therefore $ f $ shall contain an expression of the
 form $-\lambda  (Y_4 -Y_1)^{a_1}(Y_4 - Y_2)^{b_1}(Y_4 -
 Y_3)^{c_1}$ where $ a + b + c = a_1 + b_1 + c_1 $. In that case
 we also get the term $-\lambda Y_4^{a_1 + b_1 + c_1} =  -\lambda
 X_4^{(a_1 + b_1 + c_1)\gamma}$ in the expansion of $ f $ which
 cancels out $\lambda
 X_4^{(a + b + c )\gamma}$. All other expressions in the expansion
 of $ f $ have the form $ \lambda Y_4^aY_1^bY_2^cY_3^d $  where $
 (b,c,d) \neq (0,0,0)$. This expression as a monomial in $ X_i $'s
  is of the form $ \lambda X_1^{a_1}X_2^{a_2}X_3^{a_3}X_4^t $
  where $ (a_1,a_2,a_3) = (b,c,d) T_3$ for  \[ T_3 = \left(\begin{array}{ccc}
 -\delta_{11} & \delta_{12} & \delta_{13} \\
  \delta_{21} & -\delta_{22} & \delta_{23} \\
  \delta_{31} & \delta_{32} & -\delta_{33}
\end{array}\right)\] Further, $ t = a\gamma + b\delta_{14}
+c\delta_{24} + d\delta_{34}$.  As determinant of $T_3$ is not
zero, $ (a_1,a_2,a_3) \neq (0,0,0) $ if and only if $(b,c,d) \neq
(0,0,0)$. Hence the result follows.
 \epf
   We now proceed as in [3], but assume $\Im $ to be the set of
   all homogeneous polynomials $ F $ of degree $p + q $ in $ K[\pi] = K[Y]^E $
   such that $ F = f_0Y_4^q +\mbox{ terms of lower degree in } Y_4
    $,
   where powers of $ X_2 $  and $ X_3 $ are $ \geq 0 $ for any $ b
   \in Supp. F $. To complete the proof of  Lemma 2.4 (ii)  we
   need to show the existence of $ F\in \Im $ such that power of $
   X_1 $ is $ \geq 0 $ for all $ b\in Supp. F $. If this is not
   true then for any $ F\in \Im $ there exists $ b =
   (b_1,b_2,b_3,b_4)\in Supp. F $ such that $ -b_1\delta_{11} +
   b_2\delta_{21} + b_3\delta_{31}< 0  $. Then $  -b_1\delta_{11}
   +(b_2 + b_3) min(\delta_{21},\delta_{31}) < 0 $. Now, let for any
   $ F\in \Im $
    $$ e = max \{ b_4| b = (b_1,b_2,b_3,b_4)\in Supp. F,\mbox{ and }
    -b_1\delta_{11}+ (b_2 + b_3) min(\delta_{21},\delta_{31}) < 0
    \}$$ and $$ d = max \{ b_1| b = (b_1,b_2,b_3,e)\in Supp.
    F\}.$$Define as in [3], $O(F)= (d,e) $ and consider the maximal
    element $ F\in \Im $ with respect to the lexicographic
    ordering in $ \Z^2 $. Let $ h\in K[Y_2, Y_3 ] $ be the
    coefficient of $ Y_1^dY_4^e $ in $ F $. We, now, need to prove
     [3, Lemma 3.2 ]. As, in its proof in  [3], if $ \lambda
    Y_2^{c_2} Y_3^{c_3}\in E(h) $ where $ \lambda (\neq 0)\in K $,
    then as $ E(F) = 0 $ , there exists $ (d,c_2,c_3,e+1)\in Supp.
    F $. By definition of $ e $ , \[-d\delta_{11}+(c_2+c_3)\,
    min(\delta_{21},\delta_{31})\geq 0\;\;\;\cdots (i)\] We claim
    that if $(d,l_2,l_3,e)\in supp\,F,$ then $-d\delta_{11} +
    (l_2+l_3)\,min(\delta_{21},\delta_{31})<0.$ By definition of
    $e$, there exists $(a_1,a_2,a_3,e)$ in $supp\, F$ such that \[
    \begin{array}{clc}
       & -a_1\delta_{11} +
    (a_2+a_3)\,min(\delta_{21},\delta_{31})<0 &  \\
      \Rightarrow & -d\delta_{11} +
    (l_2+l_3)\,min(\delta_{21},\delta_{31})<0 & \;\;\;\cdots (ii) \\
\end{array}\] since by choice of $d$, $a_1\leq d$, and
$l_2+l_3\leq a_2+a_3.$ Hence the claim follows. From this, as
$c_2+c_3+1=l_2+l_3,$ we get \[ -d\delta_{11} +
    (c_2+c_3)\,min(\delta_{21},\delta_{31})<0.\] This contradicts
    $(i)$. Hence the proof of [3 Lemma 3.2] follows.

    For the final part of the proof, for $s=p+q-d-e$ and
    $d=\overline{p}_2+\overline{p}_3,$ consider\[G =(Y_3-Y_2)^s
    (Y_3-Y_1)^{\overline{p}_2} (Y_2-Y_1)^{\overline{p}_3}
    (Y_4-Y_1)^e\] Then $G$ is homogeneous of degree $p+q$ and
    $\lambda(Y_3-Y_2)^sY_1^dY_4^e$  appears in the expansion of $G$
    with $\lambda(\neq 0)\in K.$ To complete the argument as in [3],
    we need to show that powers of $X_2$ and $X_3$ are $\geq 0$.
    in $G$. As powers of $X_2,X_3$ in $Y_4$ and $Y_1$ are $\geq
    0$, clearly the same holds for $ (Y_4-Y_1)^e$ . Thus to prove
    our assertion it suffices to show that the same holds for
    $(Y_3-Y_2)^s(Y_3-Y_1)^{\overline{p}_2}
    (Y_2-Y_1)^{\overline{p}_3}.$ Let \[\overline{p} = s +
    \overline{p}_2+\overline{p}_3 = p+q-e =s+d.\]  Then as $e<q,$
    $\overline{p}>p.$ Therefore \[\overline{p} (1- \xi_1 -\xi_2-
    \xi_3)> p(1- \xi_1 -\xi_2-
    \xi_3)\geq 3. \;\;\;\cdots (iii)\] From the equation $(ii),$
    it is clear that \[\begin{array}{clc}
       & s\,min(\delta_{21},\delta_{31})< d\delta_{11} &  \\
      \Rightarrow & s< (s+d)\xi_1= \overline{p}\xi_1 & \;\;\;\cdots (iv) \
    \end{array}\] By the equation $(iii)$, \[ \begin{array}{clc}
       & \overline{p}-\overline{p}\xi_1-\overline{p} (\xi_2+\xi_3)\geq 3 &  \\
      \Rightarrow & s+d-\overline{p}\xi_1-\overline{p} (\xi_2+\xi_3)\geq 3 &  \\
      \Rightarrow & d-\overline{p}(\xi_2+\xi_3)\geq 3 & \mbox{  by }(iv) \
    \end{array}\] Therefore we can write $d=\overline{p}_2 +
    \overline{p}_3$ such that $\overline{p}_2> \overline{p}\xi_2$,
    $\overline{p}_3>\overline{p}\xi_3.$ Then the claim follows
    from the Lemma 3.1. The rest of the argument is same as in
    [3]. This completes the proof of our claim. $\square$

    \section*{Acknowledgement} I am thankful to Prakash Bhagat,
    discussions  with whom clarified parts of calculations in [3].

\end{document}